\documentclass[sigconf]{acmart}
\AtBeginDocument{%
  }

\setcopyright{acmlicensed}
\copyrightyear{2026}
\acmYear{2026}
\acmDOI{XXXXXXX.XXXXXXX}
\acmConference[ICAIF]{ACM INTERNATIONAL CONFERENCE ON AI IN FINANCE}{2026}{Milan, Italy}
\acmSubmissionID{215}

\usepackage{subcaption}

\begin{document}

\title{Simulation-Based Neural Policies for Portfolio Choice: Architecture, Training, and Interpretability}

\author{Jules Viard}
\correspondingauthor
\authornotemark[1]
\affiliation{%
  \institution{Imperial College London}
  \country{United Kingdom}
}
\email{viard.jules@live.fr}

\author{Alexander Michaelides}
\affiliation{%
  \institution{Imperial College London}
  \country{United Kingdom}
}
\email{a.michaelides@imperial.ac.uk}

\author{Panos Parpas}
\affiliation{%
  \institution{Imperial College London}
  \country{United Kingdom}
}
\email{p.parpas@imperial.ac.uk}

\renewcommand{\shortauthors}{Viard et al.}


\begin{abstract}
Many economic decision problems, lifecycle consumption–saving and dynamic portfolio choice, are finite-horizon stochastic control problems with continuous states and actions. When the state is low-dimensional these problems are solved by dynamic programming on a grid. The grid cost grows exponentially in the state dimension, known as the curse of dimensionality, which motivates replacing the value-function grid with a neural policy optimized directly through simulation. Such policies are usually studied in the high-dimensional settings that motivate them, precisely where no reference solution exists. So the contribution of any single architectural or training choice cannot be isolated and diagnosed. We therefore take a step back and treat both the architecture and the solution method as the objects of study. To this end, we consider a lifecycle problem with a sufficiently low-dimensional normalized state space to admit an accurate dynamic programming solution, which is used for evaluation. We compare four architectures. The simplest consists of a single time-conditioned network. We then consider two networks concatenated across the regime switch, followed by one network per date trained backward against frozen downstream policies. Finally, we evaluate a constrained variant of the per-date architecture. Decoupling the policy across time gives each date a short, well-posed objective, which we pair with direction-dominant optimization that normalizes away gradient magnitude. Architectures that lead to similar realized utility objective can nevertheless differ in whether they respect the underlying problem's economics. We therefore evaluate each design jointly based on welfare, a solution-free Bellman residual, shape restrictions, and the resulting policy functions. 
\end{abstract}

\maketitle

\section{Introduction}

Economic decision-making is, fundamentally, dynamic optimization under uncertainty: an agent chooses actions over time to maximize an expected discounted payoff subject to a stochastic law of motion using control variables~\cite{sutton2018reinforcement,rao2023foundations}.

A standard benchmark in dynamic economics is the stationary, infinite-horizon control problem: when preferences, constraints, and shocks are time-invariant, a single value function solves a Bellman equation and the optimal rule is stationary. However, lifecycle problems are time-dependent and may involve regime switches, which break the stationarity assumption. The economic environment changes over the life cycle, so the value function and the optimal policy must depend explicitly on age. 
The lifecycle portfolio problem we study is non-stationary: it has a hard terminal date, age-dependent income, and an abrupt structural break at retirement, so no time-invariant rule exists. 

Using Bellman's principle of optimality~\cite{bellman1957dynamic}, the problem can be solved recursively. The terminal period admits a closed-form solution, after which the optimal decision at each preceding date is obtained by solving a one-period optimization problem conditional on the already-characterized continuation value (value function at time $t+1$). 
\begin{equation}
    V_t(x_t) = \max_{a_t \in \mathcal{A}} \{ u(x_t, a_t) + \beta \mathbb{E}[V_{t+1}(x_{t+1}) \mid x_t, a_t] \}
\end{equation}

When the state is low-dimensional this recursion is solved on a grid with an interpolated value function. The grid cost grows exponentially in the number of state variables, so problems with several decision or state variables quickly become intractable. This curse of dimensionality is the central obstacle that motivates replacing the value-function grid with function approximation, and neural networks in particular.

We study the canonical normalized lifecycle consumption model~\cite{cocco2005consumption}. A household chooses how much to consume and what share of its savings to hold in a risky asset. The problem's defining feature is a change in the data-generating process partway through the horizon: at retirement, stochastic permanent-income growth and labor income give way to a deterministic pension.

This break motivates our central architectural question: should a single policy network represent both regimes simultaneously, or should the solution be constructed by concatenating separately trained sub-policies, each specialized to one regime? In the limit, should one use a neural network to approximate each conditional distribution at every time $t$? Because the normalized state is low-dimensional, a dynamic programming approach remains computationally feasible. We consider this solution to be the reference benchmark (degree of approximation determined by the resolution of the grid).

Our economic model is known, smooth, and continuous-action. Using the problem's structure, we optimize a deterministic policy by differentiating the realized utility along each simulated path. Given that the problem consists of a 46 period working life followed by a 35 period retirement phase, the model involves a total of 81 wealth transitions. Backpropagating the objective through the 80 transitions preceding the terminal condition weakens the direct, age-specific learning signal available for decisions made later in life. Therefore, we study different decoupling architectures, paired with training techniques for variance reduction and exploration.

\begin{itemize}
    \item \textbf{A --- single time-conditioned network}: parameters shared across all dates, with age supplied as a feature.
    \item \textbf{B --- two-regime concatenation}: a working-life network plus a separately trained, then frozen, retirement network that serves as its continuation.
    \item \textbf{C --- full backward per-date networks}: one neural network per date, trained backward, whose continuation is the realized utility of rolling out the frozen downstream policy networks.
    \item \textbf{D --- shape-constrained backward networks}: architecture C augmented with an economically motivated penalty that keeps the marginal propensity to consume positive and below one.
\end{itemize}

\paragraph{Contributions.}
\begin{itemize}
    \item A \textbf{regime-decoupled neural policy}: separate networks for working life and retirement, with the retirement policy trained first and frozen to provide the continuation value across the regime change.
    \item A \textbf{backward-induction architecture}: one network per date trained backward against frozen downstream policies, reducing each long-horizon optimization to a one-step problem.
    \item \textbf{Economically motivated constraints}: a shape penalty enforcing the theoretically implied marginal propensity to consume, adding structure information to the neural network solution.
    \item A \textbf{scalable training recipe}: payoff normalization, direction-dominant optimizers, antithetic sampling, and feasibility-preserving policy maps for high-dimensional settings.
\end{itemize}

\paragraph{Findings.}
Welfare loss, policy accuracy, and the Bellman residual improve as backward-induction structure is added from A to B to C. The single network already performs well, but it tends to over-save on average during the retirement phase. Full backward induction lands within $0.13\%$ of the certainty-equivalent reference, cuts consumption share error, and significantly reduces the fraction of paths below DP realized utility. However, the decoupling produces states in which the marginal propensity to consume become negative. The shape-constrained architecture removes them entirely, by imposing the positivity of the marginal propensity to consume that the economics of the problem dictates. For small training samples, the shape constraint also improves sample efficiency. It keeps MPC violations near zero and modestly reduces excess dispersion in realized utility.

\section{Relevant Literature}
 
The theoretical foundations of portfolio choice rest on two classical benchmarks. Markowitz’s mean-variance analysis formulates a single-period allocation problem as a trade-off between expected return and return variance, yielding the static efficient frontier but abstracting from dynamic rebalancing \cite{markowitz1952portfolio}. Merton recast portfolio choice as a continuous-time stochastic-control problem. Under CRRA preferences, constant investment opportunities, and risky returns following geometric Brownian motion, the optimal risky portfolio share is constant \cite{merton1969lifetime}. In a finite horizon setting, the consumption policy is generally time dependent.
 
We adopt the lifecycle portfolio optimization formulation of Cocco, Gomes \& Maenhout~\cite{cocco2005consumption} (surveyed in~\cite{gomes2020portfolio}): a household that chooses consumption and a risky share over a finite horizon, facing permanent and transitory labor-income shocks in working life and a deterministic pension in retirement. We use the normalization trick: CRRA homogeneity divides permanent income out of the level problem, collapsing the state to normalized cash-on-hand and lowering its dimension~\cite{carroll2022solution}. The dynamic-programming solution on a grid can be used because the normalized state is low-dimensional. We retain this solution as a reference against which every neural policy is measured.

Duarte et al.\ bring simulation-trained neural networks to high-dimensional lifecycle allocation problems~\cite{duarte2021simple}, showing that a directly parameterized policy can approximate optimal behaviour where grid-based dynamic programming is infeasible. Their setting also exposes the central methodological difficulty: once the state space is large, no reference solution exists, the policy depends on many interacting hyperparameters, and the complexity of the problem makes it hard to attribute performance to any single architectural choice or training recipe. With this controlled environment, we aim to characterize the accuracy of the neural network policy.
 
The move to neural approximation is motivated by the curse of dimensionality~\cite{bellman1957dynamic,keogh2017curse}: the grid cost of dynamic programming grows exponentially in the state dimension. Neural networks are flexible, scalable approximators (universal in the limit) ~\cite{hornik1989multilayer,goodfellow2016deep}. A growing literature applies deep learning directly to economic dynamic models. Maliar \& Winant cast lifetime-reward, Bellman, and Euler-equation objects as regression problems solved by stochastic gradient descent on simulated paths~\cite{maliar2021deep}. In the stochastic-control tradition, Han \& E approximate the control at each date by a feed-forward network and optimize the policy through Monte-Carlo simulation of the known dynamics~\cite{hanE2016deep}. 
 
The same problem is studied under the reinforcement-learning lens, which formalizes sequential decision-making as a Markov decision process~\cite{sutton2018reinforcement,rao2023foundations}. Deep reinforcement learning replaces the value or policy functions with a network, an approach behind landmark results in games~\cite{mnih2013playing,silver2016mastering}. Two axes of this literature directly inform on our design. Along the model-free/model-based axis, model-free algorithms learn from sampled reward alone. In our setting, where the model is known and smooth, we can use the pathwise gradient of the realized utility. Policy-gradient methods directly optimize a parameterized policy, ranging from the REINFORCE algorithm~\cite{williams1992simple} to deterministic policy gradients for continuous action spaces~\cite{silver2014deterministic} and their deep, off-policy implementation in DDPG~\cite{lillicrap2016continuous}. Reinforcement learning has been applied widely across finance, economics, portfolio management and trading~\cite{moody2001learning,jiang2017deep,li2019optimistic,tsang2020deep, acero2024deep}. It has also been explored in a consumption-saving setting~\cite{shi2021learning}. 
 
Our contributions lie at the intersection of dynamic portfolio choice and reinforcement learning using function approximation. We study the same class of simulation-trained neural policies as Duarte et al.~\cite{duarte2021simple}, but push the construction to a more granular, per-date full-backward architecture. This class of neural network policies is evaluated in a controlled, low-dimensional lifecycle problem. The controlled setting allows us to vary one architectural or training choice at a time, including regime decoupling, per-date backward induction, shape constraints, and training techniques. Each architectural choice can then be assessed against the reference solution obtained through dynamic programming.

\section{Problem Formulation}

\subsection{Finite-horizon stochastic control}

The lifecycle problem belongs to the class of finite-horizon stochastic control problems~\cite{hure2021deep,peng2024machine}
\begin{equation}
\max_{\pi}\ \mathbb{E}_0\!\left[\sum_{t=0}^{T} u(s_t,a_t)\right]
\end{equation}
\begin{equation*}
    a_t=\pi(s_t),\quad
s_{t+1}=m(s_t,a_t,\epsilon_{t+1}),\quad
s_0\sim F,
\end{equation*}
with a Markov state $s_t$, a continuous control $a_t$, a known, forward-simulable transition $m$ driven by exogenous shocks $\epsilon_{t+1}$, and an initial-state distribution $F$. The recursive (Bellman) view of this problem is what dynamic programming solves. The equivalent trajectory form above is what the simulation-based method of the next section optimizes directly. 

\subsection{The normalized lifecycle consumption--portfolio problem}

A household lives from age $t_b$ to $t_d$, indexing time by $t=0,\dots,T$ with $T=t_d-t_b$. Each period it chooses consumption $c_t$ and the share $\alpha_t$ of its savings held in a risky asset. Permanent income $P_t$ is divided out of the level problem by CRRA homogeneity (Section~\ref{subsec:normalization}), leaving a single continuous state: normalized cash-on-hand $x_t=X_t/P_t$. The decision state and controls are
\begin{equation}
s_t=(t,x_t),\qquad a_t=(c_t,\alpha_t),\qquad 0\le c_t\le x_t,\qquad 0\le \alpha_t\le 1 .
\end{equation}
Writing savings as $x_t-c_t$, the gross return on savings is $\mathcal{R}_{t+1}(\alpha_t)=(1-\alpha_t)\,r+\alpha_t R_{t+1}$, where the risky gross return is $R_{t+1}=r+\mu+\sigma_r\,\xi_{t+1}$ with $\xi_{t+1}\sim\mathcal{N}(0,1)$.

The problem's defining feature is a change in the data-generating process at the retirement date $t_R$. In working life ($t<t_R$) permanent-income growth $G_{t+1}$ and transitory income $Y_{t+1}$ are stochastic. At retirement ($t\ge t_R$) growth is shut off, $G_{t+1}=1$, and labor income is replaced by a fixed pension $\bar y^R$. The normalized state therefore evolves as
\begin{equation}
x_{t+1}=
\begin{cases}
\dfrac{\mathcal{R}_{t+1}(\alpha_t)\,(x_t-c_t)}{G_{t+1}}+Y_{t+1}, & t<t_R \quad\text{(working)},\\[2.2ex]
\mathcal{R}_{t+1}(\alpha_t)\,(x_t-c_t)+\bar y^R, & t\ge t_R \quad\text{(retirement)}.
\end{cases}
\end{equation}

\subsection{Normalization and the payoff weight}\label{subsec:normalization}

Period utility is constant relative risk aversion, $u(c)=c^{1-\rho}/(1-\rho)$ with $\rho\neq1$. Its homogeneity of degree $1-\rho$ lets permanent income be divided out~\cite{carroll2022solution}, collapsing the problem to the one-dimensional state $x_t$ but leaving a residual continuation factor $G_{t+1}^{1-\rho}$ in working life. The normalized Bellman equation then reads
\begin{equation}
v_t(x)=\!\!\max_{\substack{0\le c\le x\\ 0\le\alpha\le1}}\!\Big\{\tfrac{c^{1-\rho}}{1-\rho}+\delta\,\mathbb{E}_t\big[\kappa_{t+1}\,v_{t+1}(x_{t+1})\big]\Big\}
\end{equation}
\begin{equation*}
    \kappa_{t+1}=\begin{cases}G_{t+1}^{1-\rho}, & t<t_R,\\[0.4ex] 1, & t\ge t_R,\end{cases}
\end{equation*}
with $x_{t+1}$ given by the regime-appropriate transition and the terminal condition $v_T(x)=x^{1-\rho}/(1-\rho)$ (e.g. full consumption).

Along a simulated path it is accumulated into a payoff weight $d_t$, with $d_0=1$ and $d_{t+1}=d_t\,\delta\,\kappa_{t+1}$, so that maximizing the normalized lifetime objective
\begin{equation}
\mathbb{E}\!\left[\sum_{t=0}^{T} d_t\,\frac{c_t^{1-\rho}}{1-\rho}\right]
\label{eq:objective}
\end{equation}
is equivalent to solving the Bellman equation above. 

\begin{table}[t]
\centering\small
\caption{Model calibration.}
\label{tab:calibration}
\begin{tabular}{llr}
\toprule
Component & Symbol & Value\\
\midrule
Relative risk aversion & $\rho$ & $5.0$ \\
Discount factor & $\delta$ & $0.97$ \\
Gross risk-free return & $r$ & $1.015$ \\
Mean excess return & $\mu$ & $0.04$ \\
Return volatility & $\sigma_r$ & $0.20$ \\
Shock correlations & \multicolumn{2}{r}{$\mathrm{corr}(\varepsilon,\nu)=\mathrm{corr}(\cdot,\xi)=0$}\\
Retirement income & $\bar y^R$ & $0.68212$ \\
Initial normalized cash & $x_0$ & $Y_0$ (no initial wealth) \\
State Cash grid & $x$ & $[0.25,115]$ \\
\bottomrule
\end{tabular}
\end{table}

Table~\ref{tab:calibration} collects the economic calibration used throughout the benchmark and neural-policy experiments.

\section{Simulation-Based Neural Policy Optimization}

\subsection{Differentiable direct policy optimization}

We parameterize the policy as a deterministic neural network $a_t=\pi_\theta(t,x_t)$ and optimize it by forward simulation, as in~\cite{duarte2021simple}. The simulated objective is a deterministic, differentiable function of $\theta$,
\begin{equation}
\widehat J(\theta)=\frac{1}{N}\sum_{i=1}^{N}\sum_{t=0}^{T} d_t^{(i)}\,\frac{\big(c_t^{(i)}(\theta)\big)^{1-\rho}}{1-\rho},
\label{eq:mc-objective}
\end{equation}
the sample analogue of the normalized objective in Eq.~\eqref{eq:objective}. Each path is rolled forward under $\pi_\theta$ through the working/retirement transition introduced above, and $\theta$ is updated by stochastic gradient ascent on $\widehat J$ over minibatches of paths.

The gradient is the pathwise estimator, obtained by automatic differentiation of the rollout,
\begin{equation}
\nabla_\theta \widehat J(\theta)=\frac{1}{N}\sum_{i=1}^{N}\sum_{t=0}^{T} d_t^{(i)}\,\big(c_t^{(i)}\big)^{-\rho}\,\frac{\mathrm{d} c_t^{(i)}}{\mathrm{d}\theta}.
\label{eq:pathwise-grad}
\end{equation}
The realized cash $x_t^{(i)}(\theta)$ is itself a function of $\theta$ through every earlier decision, so the date-$t$ term propagates through the entire trajectory.

\subsection{Feasibility by construction}

The constraints $0\le c_t\le x_t$ and $0\le\alpha_t\le1$ are enforced by the policy map. In our trained models the consumption share is $c_t=\big(\phi+(1-\phi)\,\sigma(z_1)\big)\,x_t$ for a small floor $\phi$ and logistic $\sigma$, and the risky share is $\alpha_t=\sigma(z_2)$. The terminal rule $c_T=x_T,\ \alpha_T=0$ is imposed. Consequently the optimizer faces an unconstrained problem in $\theta$.

\subsection{Variance reduction and direction-dominant updates}\label{subsec:optim}

Two main tricks stabilize the pathwise gradient. First, the exogenous shocks are drawn once as antithetic pairs: every base path is mirrored by negating its Gaussian innovations and the pair is averaged, which lowers gradient variance. Second, the per-date objectives of the decoupled architectures differ by orders of magnitude in scale. Late-life payoff weights $d_t$ are tiny relative to early-life ones. We therefore adopt a direction-dominant update strategy. The stage update is normalized prior to the optimizer step, preserving its direction while allowing the learning rate, rather than the payoff scale, to determine the effective step size. The normalized update is then applied using AdamW. The single and two regime networks, whose objective uses multiple dates (full horizon or regime horizon), use Adamax directly.

\section{Policy Architectures}

\subsection{Architecture A --- single time-conditioned network}

\paragraph{Construction.}
A single multilayer perceptron $\pi_\theta:(t,x)\mapsto(c,\alpha)$ represents the policy at every date. One parameter vector $\theta$ is shared across the whole horizon and age enters as a feature. The network maps a compact three-feature encoding of the decision state,
\begin{equation}
\zeta(t,x)=\Big(\,t/T,\ \tfrac{\log(1+x)}{\log(1+\bar x)},\ \mathbf{1}[t\ge t_R]\,\Big),\qquad \bar x=115,
\end{equation}
through two hidden layers of $128$ $\tanh$ units to two raw outputs, which the feasibility map (floor $\phi=0.005$) turns into $(c,\alpha)$. It has $17{,}282$ trainable parameters.

\paragraph{Training objective and algorithm.}
The parameters $\theta$ are trained to maximize the full-horizon simulated objective $\widehat J(\theta)$ of Eq.~\eqref{eq:mc-objective} directly. A fixed set of antithetic shock paths is drawn once; each is rolled forward under $\pi_\theta$ across the working/retirement transition. The payoff-weighted lifetime utility is averaged, and $\theta$ is updated by stochastic gradient ascent on minibatches using the pathwise gradient of Eq.~\eqref{eq:pathwise-grad}. The reported model uses $100{,}000$ base paths (doubled by antithetic mirroring), minibatches of $2{,}048$, $2{,}000$ epochs, and the Adamax optimizer at learning rate $5\times10^{-3}$.

\paragraph{Deployment.}
At any decision date the action is one forward pass $\pi_\theta(t,x_t)$. The same network serves both accumulation and decumulation of the wealth. The regime indicator being its only explicit signal of the retirement break.

\subsection{Architecture B --- two-regime concatenation}

\paragraph{Construction.}
The change in the underlying law of motion at retirement motivates representing the two regimes by two networks rather than one smooth function of age. A working-life network $\pi_{\theta_W}$ governs the dates $t<t_R$ and a retirement network $\pi_{\theta_R}$ governs $t\ge t_R$. Each uses the same three-feature encoding $\zeta$ and the same depth as Architecture A but with $64$-unit hidden layers, giving $4{,}546$ parameters per network ($9{,}092$ in total).

\paragraph{Training objective and algorithm.}
Training is backward in two blocks: \emph{(i) Retirement}. $\pi_{\theta_R}$ is trained on retirement only paths that begin from a sampled retirement-entry cash $x_{t_R}$. This policy should maximize the retirement-block objective
\begin{equation}
\widehat J_R(\theta_R)=\frac1N\sum_{i=1}^N\sum_{t=t_R}^{T}\tilde d_t^{(i)}\,\frac{\big(c_t^{(i)}\big)^{1-\rho}}{1-\rho},
\end{equation}
The reported model draws $x_{t_R}$ entirely from a broad log-uniform distribution over the cash grid, so the retirement network sees no DP information. \emph{(ii) Working life.} $\theta_R$ is then frozen and $\pi_{\theta_W}$ is trained on full lifecycle paths in which every action from $t_R$ onward is produced by the frozen $\pi_{\theta_R}$. Thus, $\theta_W$ maximizes the full objective $\widehat J(\theta)$ of Eq.~\eqref{eq:mc-objective} against a fixed continuation environment. Both networks are optimized with Adamax at learning rate $5\times10^{-3}$ for $4{,}000$ epochs over $100{,}000$ base paths (antithetic-mirrored) and a batch size of $2{,}048$.

\paragraph{Deployment.}
The combined policy dispatches by age, $\pi(t,x)=\pi_{\theta_W}(t,x)$ for $t<t_R$ and $\pi_{\theta_R}(t,x)$ for $t\ge t_R$ (here $t_R$ is date $46$). It represents the two-block analogue of backward induction.

\subsection{Architecture C --- full backward per-date neural networks}

\paragraph{Construction.}
Backward-induction structure is carried to every date. Each date $t$ owns its own policy network $\pi_{\theta_t}$. The parameter vectors $\theta_0,\dots,\theta_T$ are disjoint, so time is represented by the network index. Each neural network maps the single cash feature $\log(1+x)/\log(1+\bar x)$ through two hidden layers of $32$ $\tanh$ units to the feasibility map (floor $\phi=0.005$). Each neural network has $1{,}186$ parameters. The terminal network is determined by the final condition: $\pi_{\theta_T}(x)=(x,0)$.

\paragraph{Training objective and algorithm.}
The neural networks are trained following a backward methodology: $t=T-1,\dots,0$. At stage $t$ all downstream networks are frozen and the continuation is the realized utility of rolling those frozen networks forward to $T$
\begin{equation}
\widehat J_t(\theta_t)=\frac1N\sum_{i=1}^N\sum_{s=t}^{T}\tilde d_s^{(i)}\,\frac{\big(c_s^{(i)}\big)^{1-\rho}}{1-\rho},
\end{equation}
in which the date-$t$ action comes from the trainable neural network and all later actions from the frozen policy networks. The start cash $x_t^{(i)}$ is drawn from a broad log-uniform distribution over the cash grid. Once converged, $\theta_t$ is frozen and the training moves to $t-1$. Each neural network is trained for $3{,}500$ epochs over $200{,}000$ base paths (antithetic-mirrored) in minibatches of $512$, using the direction-dominant update of Section~\ref{subsec:optim} (the stage gradient normalized to unit norm, then an AdamW step at learning rate $10^{-3}$).

\paragraph{Deployment.}
The trained policy is the frozen sequence $\{\pi_{\theta_t}\}_{t=0}^{T}$. At a decision date it dispatches to the policy network whose index matches the integer model time, with the analytic network at $T$ (terminal condition).

\subsection{Architecture D --- shape-constrained backward neural networks}

\paragraph{Motivation.}
Economic theory pins down the shape of the consumption rule even where no reference solution is available: consumption is nondecreasing in cash-on-hand and rises by at most one unit per unit of cash. The marginal propensity to consume (MPC) satisfies $0\le \partial c/\partial x\le 1$~\cite{carrollkimball1996concavity, carroll2022solution}. Because each per-date objective is nearly flat in the local curvature of the consumption rule, many nearby rules deliver almost the same one-step value. The unconstrained policy networks leave the MPC weakly identified and, as Section~\ref{sec:results} reports, generate the largest number of states with a negative MPC among all four designs.

\paragraph{Construction and training.}
Architecture D follows the same structure and parameterization as Architecture C. It adds, to each stage objective, a penalty for MPC violations. At stage $t$ the MPC $\partial c/\partial x$ is obtained by automatic differentiation of the network's consumption at the sampled cash points. The stage loss becomes
\begin{equation}
-\,\widehat J_t(\theta_t)\;+\;\lambda\,\mathbb{E}\big[\,\mathrm{relu}(-\partial c/\partial x)^2+\mathrm{relu}(\partial c/\partial x-1)^2\,\big],
\end{equation}
with a weight $\lambda$ that makes positivity and the upper bound effectively binding. Enforcing the constraint by a penalty rather than by construction keeps the network, feasibility map, sampler, and optimizer identical to Architecture C.

\paragraph{Deployment.}
Deployment is identical to Architecture C.

\subsection{Architecture comparison}

Table~\ref{tab:architectures} summarizes the four designs. They differ in how time is represented. Age is an input feature in A, a dispatch between two networks in B, and the policy network index in C and D. Decoupling trades parameter sharing for well-posed per-date objectives. Architectures C and D differ because of the MPC penalty. That isolates the effect of the shape constraint.

The different architectures represent distinct model types and solution methods and require independent optimization and hyperparameter tuning. Each architecture was tuned separately, and the reported results correspond to the best configuration we found for that architecture under a comparable computational budget. The comparison therefore evaluates well-tuned representatives of each design, allowing us to distinguish the benefits of the different solution methods. Since Architectures C and D belong to the same class of model and use the same solution method, they share the same hyperparameter settings. This allows us to isolate the impact of the MPC constraint.

The replication code will be made available on GitHub, including the full calibration, grid, training, and evaluation parameters used for each reported architecture.

\begin{table*}[t]
\centering\small
\caption{Policy architectures. All four are trained without reference to the DP solution.}
\label{tab:architectures}
\begin{tabular}{lllll}
\toprule
 & \textbf{A --- Single} & \textbf{B --- Two-regime} & \textbf{C --- Full backward} & \textbf{D --- Shape-constrained}\\
\midrule
Policy objects & $1$ network & $2$ networks & $81$ networks ($80{+}$ term.) & $81$ networks ($80{+}$ term.)\\
Parameters & $17{,}282$ & $9{,}092$ ($2\times4{,}546$) & $94{,}880$ ($80\times1{,}186$) & $94{,}880$ ($80\times1{,}186$)\\
Hidden layers & $2\times128$ & $2\times64$ & $2\times32$ & $2\times32$\\
Input features & $3$ & $3$ & $1$ (cash only) & $1$ (cash only)\\
Time representation & network feature & regime dispatch & network index & network index\\
Continuation & implicit rollout & frozen retirement net & frozen downstream networks & frozen downstream networks\\
Active optimizer & Adamax & Adamax & norm.\ update $\to$ AdamW & norm.\ update $\to$ AdamW\\
Entry-cash sampling & endogenous & broad (retirement) & broad (all stages) & broad (all stages)\\
Shape constraint & none & none & none & MPC $\in[0,1]$ penalty\\
\bottomrule
\end{tabular}
\end{table*}

\section{Results}\label{sec:results}

\subsection{Welfare and pathwise accuracy}

Welfare is reported in Table~\ref{tab:welfare}. All four policies are near-optimal and dominate naive rules by two orders of magnitude (consume-all loses $12.7\%$ of CE consumption, equal-allocation $44.9\%$). Certainty-equivalent loss falls monotonically as backward-induction structure is added: the two-regime split closes about $40\%$ of the single network's CE loss. The full backward induction closes just over half, landing within $0.13\%$ CE consumption of the reference solution. The shape penalty of D does not measurably change it. Under A and B roughly four paths in five fall below the DP solution. The full-backward models cut this to $58\%$ (C) and $54\%$ (D) and shrink the median path gap. The deepest-loss paths (fifth-percentile gap) are comparable across designs, mildest for the two-regime net.

\begin{table}[t]
\centering\small
\caption{Welfare and pathwise accuracy versus the DP reference. CE loss is in certainty-equivalent consumption units; gaps are NN${-}$DP. Bold marks the best of the four.}
\label{tab:welfare}
\begin{tabular}{lcccc}
\toprule
 & \textbf{A} & \textbf{B} & \textbf{C} & \textbf{D}\\
 & Single & Two-reg. & Full bw. & Shape\\
\midrule
CE loss vs.\ DP        & $-0.269\%$ & $-0.167\%$ & $\mathbf{-0.127\%}$ & $-0.127\%$\\
Objective gap vs.\ DP  & $-0.0278$  & $-0.0173$  & $\mathbf{-0.0130}$   & $-0.0131$\\
\midrule
Paths below DP         & $79.7\%$   & $79.1\%$   & $58.4\%$   & $\mathbf{53.8\%}$\\
Median path gap        & $-0.019$   & $-0.009$   & $-0.004$   & $\mathbf{-0.002}$\\
5th-pct.\ path gap     & $-0.110$   & $\mathbf{-0.069}$   & $-0.090$   & $-0.095$\\
\bottomrule
\end{tabular}
\end{table}

\subsection{Policy Accuracy}

Pointwise policy error against DP, on the shared cash grid, shows the share mean absolute error (MAE, Table~\ref{tab:mae}). Consumption share error decreases from A to the full-backward models ($0.055\to0.018\to0.013$). The remaining dimensions are noisier and non-monotone. The objective is nearly flat in the portfolio share, B and C attain the lowest risky-share MAE ($\approx0.037$) while D is slightly higher.

\begin{figure*}[t]
    \centering

    \begin{subfigure}{0.8\textwidth}
        \centering
        \includegraphics[width=\linewidth]{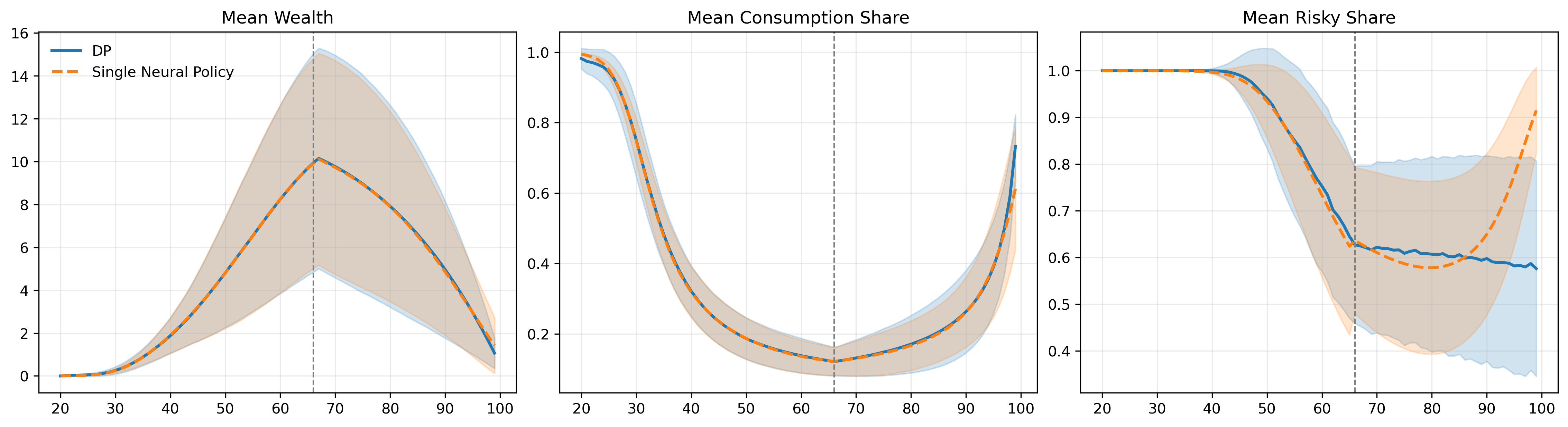}
        \caption{Policy function Architecture A}
    \end{subfigure}

    \vspace{1em}

    \begin{subfigure}{0.8\textwidth}
        \centering
        \includegraphics[width=\linewidth]{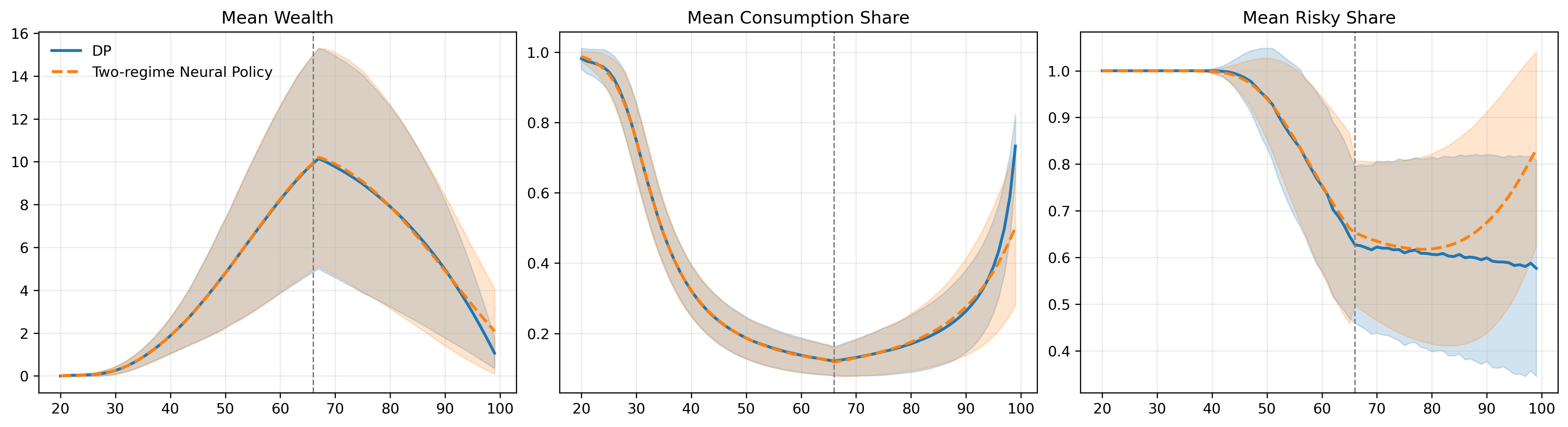}
        \caption{Policy function Architecture B}
    \end{subfigure}

    \vspace{1em}

    \begin{subfigure}{0.8\textwidth}
        \centering
        \includegraphics[width=\linewidth]{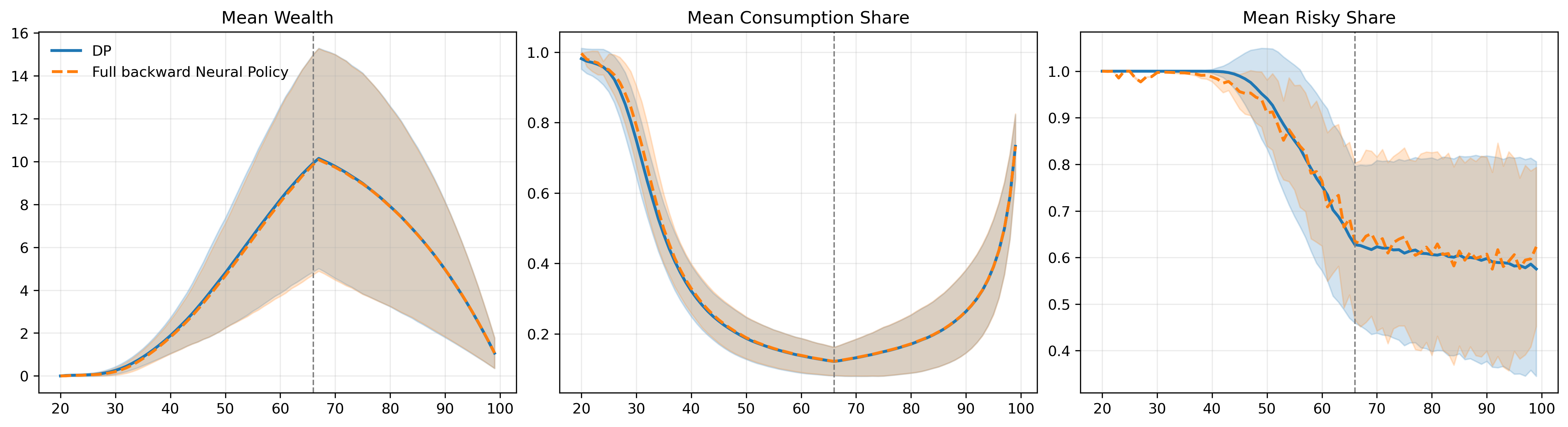}
        \caption{Policy function Architecture C}
    \end{subfigure}

    \vspace{1em}

    \begin{subfigure}{0.8\textwidth}
        \centering
        \includegraphics[width=\linewidth]{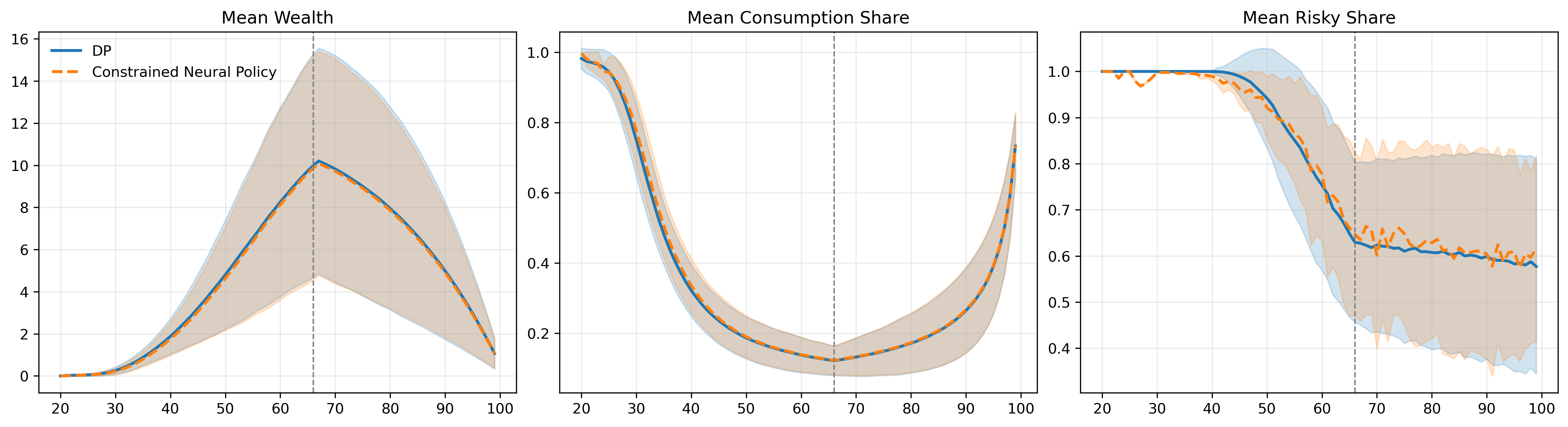}
        \caption{Policy function Architecture D}
    \end{subfigure}

    \caption{Policy panels across architectures with mean and standard deviation across 20,000 simulations.}
    \label{fig:policy-panel}
\end{figure*}

The simulated mean profiles (Fig.~\ref{fig:policy-panel}) provide an overall view of the accuracy at each date $t$. Mean normalized wealth is hump-shaped and peaks at retirement. The mean consumption share is U-shaped over the life cycle. For every architecture, the policies during the working life are well aligned with the DP solution. The different architectures affect the risky asset share during the wealth decumulation phase. The single network (A) matches DP through accumulation but drifts upward in late retirement, over-weighting the risky asset. Full backward induction (C) removes this divergence and follows the DP decline across retirement. The shape-constrained model (D) preserves that fit. The close overlap of the shaded percentile bands shows the agreement holds across the cross-sectional distribution of households.

\begin{table}[t]
\centering\small
\caption{Mean absolute policy error versus DP on the shared $(\text{age}\times\text{cash})$ grid, over nonterminal states.}
\label{tab:mae}
\begin{tabular}{lcccc}
\toprule
MAE vs.\ DP & \textbf{A} Single & \textbf{B} Two-reg.\ & \textbf{C} Full bw.\ & \textbf{D} Shape\\
\midrule
Consumption share & $0.055$ & $0.018$ & $0.013$ & $\mathbf{0.013}$\\
Risky share       & $0.069$ & $\mathbf{0.037}$ & $0.037$ & $0.047$\\
Risky savings     & $2.28$  & $\mathbf{1.22}$  & $1.30$  & $1.73$\\
\bottomrule
\end{tabular}
\end{table}

\subsection{Solution-free diagnostics}

\paragraph{The Bellman residual (optimality).}
Our primary solution-free instrument is the one-step Bellman improvement gap~\cite{sutton2018reinforcement}. For a policy $\pi$ with value $v_\pi$, the date-$t$ action value is: $Q_t(x,a)=u(x,a)+\delta\,\mathbb{E}_t[\kappa_{t+1}\,v_\pi(x')\mid x,a]$. Therefore, the residual is
\begin{equation}
\Delta_t(x)=\max_{a}Q_t(x,a)\;-\;Q_t\big(x,\pi_t(x)\big)\ \ge\ 0,
\label{eq:bellman-residual}
\end{equation}
This represents the welfare an agent would gain by deviating optimally for one period and then reverting to $\pi$. It is nonnegative by construction and zero everywhere if and only if $\pi$ is optimal. We compute it from the trained policy alone. The continuation $v_\pi$ is estimated by Monte-Carlo simulation: from each grid state we roll the frozen policy network forward to the terminal date on the common shock paths and average the realized normalized utility of Eq.~\eqref{eq:objective}. Then, that value is interpolated with a cubic spline. The maximization uses the known one-period model (Gaussian quadrature over the shocks) with a grid search over $(c,\alpha)$. We report the gap in certainty-equivalent consumption units, $(Q^\star/Q_\pi)^{1/(1-\rho)}-1$. Because a uniform average over the cash grid can be dominated by rarely-visited corners, we also report a visitation-weighted version. 

On the uniform grid, adding backward structure reduces the CE residual: mean $2.9\%$ (A) and $2.4\%$ (B) down to $0.54\%$ (for C, D). The worst-case gap went from over $400\%$ to about $120\%$. Under visitation weighting the full-backward models are the best performing architectures (mean $0.004\%$ CE). A ($0.20\%$) and especially B ($0.57\%$) carry a heavier tail. Per-date backward induction replaces a single long-horizon objective with a sequence of short, well-posed one-step problems. It eliminates the decumulation suboptimality. The shape penalty in D leaves both the visited and uniform mean residuals essentially unchanged ($0.0042$ vs $0.0045$, $0.544$ vs $0.544$). Its effect appears in the shape diagnostics rather than in the residual level.

\paragraph{Economic shape (MPC and monotonicity).}
Feasibility ($0\le c\le x$, $0\le\alpha\le1$) holds with zero violations for all four, confirming the by-construction policy map. Economic shape degrades as time is decoupled and each per-date objective becomes locally flat in the consumption level (Table~\ref{tab:diag}). The single and two-regime nets carry a modest number of negative MPC (A: $32$ cells, $0.98\%$; B: $34$ cells, $1.04\%$). The unconstrained full-backward model C has $109$ negative-MPC cells ($3.32\%$) and $97$ non-monotone steps spread across the horizon (worst $-0.34$), because per-date decoupling removes the cross-age smoothing. Architecture D targets this issue using a soft constraint. The penalty drives negative-MPC and non-monotone cells to zero ($0\%$).

\begin{table}[t]
\centering\small
\caption{Solution-free diagnostics: Bellman residuals in certainty-equivalent consumption units (percent). Under the same diagnostic the exact DP policy scores $0.00004$ (visited median), $0.0002$ (visited mean), $0.005$ (uniform mean), and $0.10$ (uniform max). Bold marks the best of the four per row.}
\label{tab:diag}
\begin{tabular}{lcccc}
\toprule
Bellman resid.\ (CE \%) & \textbf{A} & \textbf{B} & \textbf{C} & \textbf{D}\\
\midrule
\quad visited mean   & $0.202$  & $0.569$  & $\mathbf{0.0042}$ & $0.0045$\\
\quad uniform mean   & $2.89$   & $2.38$   & $0.544$  & $\mathbf{0.544}$\\
\quad uniform max    & $412$    & $608$    & $122$    & $\mathbf{117}$\\
\midrule
Negative-MPC cells   & $32$ & $34$ & $109$ & $\mathbf{0}$\\
Non-monotone steps   & $26$ & $25$ & $97$  & $\mathbf{0}$\\
\bottomrule
\end{tabular}
\end{table}

\subsection{Sample efficiency and the shape constraint}

We retrain Architectures C and D identically, using only $100$ base paths per policy network instead of $200{,}000$. Both deteriorate comparably in welfare, from $0.13\%$ to about $1.0\%$ certainty-equivalent loss. Its effect is on the shape of the policy and on the dispersion of realized outcomes: the mean Bellman residual on the uniform grid is slightly lower for D ($0.34\%$ vs $0.39\%$ CE), and the cross-sectional standard deviation of realized lifetime utility falls from $2.305$ ($+7.0\%$ excess dispersion relative to DP) to $2.286$ ($+6.1\%$). 

With data withdrawn, the unconstrained model also deteriorates in terms of negative MPC cells: from $109$ ($3.32\%$) to $151$ ($4.60\%$) negative-MPC cells and from $97$ ($2.96\%$) to $143$ ($4.36\%$) non-monotone steps (worst MPC $-0.73$). Architecture D limits this to $4$ negative-MPC cells ($0.12\%$) and a single non-monotone step (($0.03\%$)). The shape penalty allows one to incorporate structural information directly into the policy network, thereby improving sample efficiency. Such constraint integration may be especially important in high-dimensional settings, where generating a sufficiently comprehensive training sample is substantially more computationally intensive.

\subsection{Cross-architecture synthesis}

The complete backward training configuration performs better than the single-network and two-regime configurations in welfare, central pathwise accuracy, and consumption-share accuracy. It reduces the fraction of paths that fall below DP, removing the over-saving bias. 

Introducing a two-regime split recovers roughly $40\%$ of the welfare loss at low computational cost, but its single frozen retirement block leaves the decumulation tail unresolved. This issue is eliminated only by full backward induction, which trains each retirement date separately. Decoupling the problem comes at the expense of economic structure: the unconstrained full-backward model (C) exhibits the largest number of negative-MPC regions. However, the shape-constrained Architecture D removes these violations. This constraint becomes more valuable under a scarce training sample, and potentially for future high dimensional application.

All four architectures are trained without access to the dynamic programming solution. Taken together, realized welfare, the Bellman residual, and economic-shape diagnostics provide complementary perspectives. Welfare measures realized lifecycle outcomes, the residual identifies local gains from deviation, and shape diagnostics reveal economically wrong behavior that similar objective values can conceal.

All results are based on a single training run per architecture and evaluated on a common set of shock paths. The use of common random numbers substantially reduces Monte Carlo error in cross-architecture comparisons. A systematic assessment of training-run dispersion is left to future work, focused on individual model type (solving method).

\section{Conclusion}

We studied how neural-policy architecture and the associated solution method affect lifecycle consumption and portfolio choice. Although neural methods target high-dimensional problems where grid-based DP is infeasible, we use a controlled benchmark solely for evaluation and policy-function inspection.

The main methodological gain comes from adding backward induction structure. A single time conditioned network solves one long horizon problem. A two regime design freezes the retirement policy and uses it as the continuation for working life. A full backward design trains one neural network per date against frozen downstream networks. This progression improves certainty equivalent loss from $0.269\%$ to $0.167\%$ and then to $0.127\%$. It also reduces consumption share error from $0.055$ to $0.018$ and then to $0.013$. The fraction of paths below the dynamic programming benchmark falls most clearly once full backward induction is used. 

Temporal decoupling can weaken the economic shape. Although the unconstrained and shape-constrained full-backward policies have similar certainty-equivalent losses, the unconstrained policy produces the highest number of negative MPC cells. Architecture D adds a penalty for violations of $0 \le \partial c/\partial x \le 1$. This removes negative MPC and non monotone cells in the full sample. Architecture D also records the lowest fraction of paths below the reference, showing that full backward induction and the economic shape constraint improve solution reliability. Realized utility alone is insufficient to validate a policy, as economically important local behavior can remain weakly identified by the aggregate objective. The shape constraint is also useful when training paths are scarce. The constrained policy keeps MPC violations near zero and reduces excess dispersion in realized utility.

The Bellman residual provides a solution-free measure of local suboptimality, while economic-shape diagnostics and policy function plots reveal behavior that aggregate welfare can conceal. Taken together, the results support full-backward training with explicit economic constraints. More broadly, they demonstrate the value of evaluating architectures and solution methods in a controlled setting.

\bibliographystyle{ACM-Reference-Format}
\bibliography{sample-base}

\end{document}